\documentclass[12pt,reqno]{amsart}
\usepackage{amsmath,amssymb,amsfonts,amscd,latexsym,amsthm,mathrsfs}
\usepackage[usenames]{color}
\usepackage[unicode]{hyperref}
\usepackage{graphicx}
\begin{document}

\title{A BBP\protect\lowercase{-style computation for $\pi$ in base 5}}

\author{W\protect\lowercase{adim} Z\protect\lowercase{udilin}}
\address{Department of Mathematics, IMAPP, Radboud University, PO Box 9010, 6500~GL Nij\-me\-gen, Netherlands}
\urladdr{https://www.math.ru.nl/~wzudilin/}

\date{17 September 2024}


\begin{abstract}
We joke about how to compute (promptly) the digits of $\pi$, in base 5, from a given place without computing preceding ones.
\end{abstract}

\maketitle

\section{The BBP algorithm}
\label{s1}

Given a base $b\in\mathbb Z_{\ge2}$ and a `related to the base' series expansion of constant $\xi\in\mathbb R$ (a BBP-type formula), the Bailey--Borwein--Plouffe algorithm \cite{BBP97} computes $r$ base~$b$ digits of $\xi$ beginning at the position $d+1$ after the floating point.
A traditional illustrative example is the constant
\[
\log2=\sum_{n=1}^\infty\frac1{n2^n}
\]
and $b=2$.
To access what begins at the binary position $d+1$, we look for the fractional part of
\[
2^d\log2=\sum_{n=1}^d\frac{2^{d-n}}{n}+\sum_{n=d+1}^\infty\frac{2^{d-n}}{n}.
\]
In the first (finite) sum fractional parts of individual terms are computed using
\[
\frac{2^{d-n}}{n}\bmod\mathbb Z=\frac{2^{d-n}\bmod n}n,
\]
with $2^{d-n}\bmod n$ calculated with the help of fast modular exponentiation; the second sum\,---\,the tail\,---\,converges quickly at rate $1/2$, so that its first few digits can be easily computed.

On the other hand, one can vary this recipe slightly by writing
\[
2^d\log2=\sum_{n=1}^d\frac{2^{d-n}}{n}+\sum_{n=d+1}^{d+2r}\frac{2^{d-n}}{n}+\sum_{n=d+2r+1}^\infty\frac{2^{d-n}}{n}.
\]
For the first sum here we do exactly the same as before\,---\,compute individual fractional parts with the help of fast exponentiation. For the second sum we just sum up the corresponding $2r$ fractions, each being the reciprocal of an integer, while for the third sum we use the estimate
\[
\bigg|\sum_{n=d+2r+1}^\infty\frac{2^{d-n}}{n}\bigg|<\frac1{d+2r+1}\sum_{n=d+2r+1}^\infty2^{d-n}=\frac{2^{-2r}}{d+2r+1}
\]
showing that this tail will not in general affect the first $r$ binary digits of $\{2^d\log2\}$.

\section{A BBP-type formula for $\pi$}
\label{s2}

Several BBP-type formulae for $\pi$ are known for bases which are powers of 2.
This gives access, for example, to computing hexadecimal digits of $\pi$ starting at a particular far-away position without computing predecessors.
It is an open question whether this type of computation is possible for $\pi$ related to other bases.

In what follows $i$ stands for $\sqrt{-1}$.

Rational approximations to $\pi$ constructed by Salikhov in \cite{Sa08,Sa10} and later refined by Zeilberger and this author in~\cite{ZZ20} are based on the representation
\begin{equation}
\sum_{n=0}^\infty\frac1{2n+1}\bigg(\frac1{1-2i}\bigg)^{2n+1}
-\sum_{n=0}^\infty\frac1{2n+1}\bigg(\frac1{1+2i}\bigg)^{2n+1}
=\frac{\pi i}4.
\label{eq1}
\end{equation}
Because the norm of the quadratic irrationalities $1\pm 2i$ is 5, one can think of the expansion as of `base $\sqrt5$' BBP-type formula for $\pi$\,---\,in fact `base 5' as we encounter the powers of $(1\pm 2i)^2$ rather than of $1\pm 2i$.
Can it be used for computing base~5 digits of~$\pi$?

Before proceeding with this, we make some related comments about computation of powers $(1\pm2i)^n$.
Their real and imaginary parts are read off from the $2\times2$ matrix $\big(\begin{smallmatrix} 1 & -2 \\ 2 & 1 \end{smallmatrix}\big)^n$; they are also generated via the recurrence equation $a_n=2a_{n-1}-5a_{n-2}$.
The latter circumstance leads to the following interpretation of identity~\eqref{eq1}: Define the sequence $b_n$ through the recursion $b_n=-6b_{n-1}-25b_{n-2}$ for $n\ge2$ and initial data $b_0=1$, $b_1=-1$.
Then
\[
\sum_{n=0}^\infty\frac1{2n+1}\,\frac{b_n}{5^{2n}}
=\frac{5\pi}{16}
\]
and $|b_n|<2\cdot5^n$ for all $n$.

Formula \eqref{eq1} allows one to cast $\pi$ as the imaginary part of
\[
\xi=8\sum_{n=0}^\infty\frac1{2n+1}\bigg(\frac1{1-2i}\bigg)^{2n+1},
\]
and we can use the BBP strategy to compute the base~$5$ expansion for both its real and imaginary parts beginning at the position $d+1$ till the position $d+r$, say, under a mild condition we give below (see condition~\eqref{eq2}). This demands for computing first $r$ base~5 digits after the floating point of the number $5^d\xi$ (again we do this for both real and imaginary parts!).
Because $|1-2i|=\sqrt5$, the tail
\[
8\cdot5^d\sum_{n=d+2r}^\infty\frac1{2n+1}\bigg(\frac1{1-2i}\bigg)^{2n+1}
\]
is bounded above by
\[
\frac{8\cdot5^d}{2d+4r+1}\sum_{n=d+2r}^\infty\frac1{5^{n+1/2}}
=\frac{2\cdot5^{1/2-2r}}{2d+4r+1},
\]
which means that the positions before the position $r+1$ in the base $5$ expansion of $5^d\xi$ are hardly affected, and we look for the first $r$ digits of
\[
\sum_{n=0}^{d+2r-1}\frac{8\cdot5^d}{2n+1}\bigg(\frac1{1-2i}\bigg)^{2n+1}.
\]
We are only interested in this expression in $\mathbb Q[i]$ modulo~1 for real and imaginary parts.
Assuming that
\begin{equation}
\frac{\log(2d+4r-1)}{\log5}\le d
\label{eq2}
\end{equation}
and writing $2n+1=5^k\cdot m$ with $(m,5)=1$ for the $n$th term in the sum, we look for
\begin{equation}
\frac{8\cdot5^d}{2n+1}\bigg(\frac1{1-2i}\bigg)^{2n+1}\bmod\mathbb Z[i]
=\frac{8\cdot5^{d-k}(1-2i)^{-(2n+1)}\bmod m}{m};
\label{eq3}
\end{equation}
thus, we only need executing fast exponentiation for $5^{d-k}\bmod m$ and $(1-2i)^{2n+1}\bmod m$ in $\mathbb Z[i]$, and then computing the reciprocal of the latter modulo~$m$. The latter inversion is possible because $m$ is coprime with the norm of any power of $1-2i$; alternatively, one can take an integer $a$ such that $a\equiv\frac15\bmod m$, with the motive that $(1-2i)^{-1}=\frac15+\frac25i\equiv a(1+2i)\bmod{m\mathbb Z[i]}$, and compute $a^{2n+1}(1+2i)^{2n+1}\bmod m$.

\section{An obvious flaw}
\label{s3}

The equality in \eqref{eq3} is incorrect when $2n+1>d-k$, since the left-hand side in the latter case has the denominator $5^{2n+1-d+k}m$, while the denominator of the right-hand side is~$m$.
Using a shorter finite sum over $n$ to meet this constraint is not an option either, because there is no bound for (the fractional part of) the tail in such cases.
Are there more suitable formulae for $\pi$ or other interesting constants, defined over quadratic or more general algebraic extensions, that could be of use?



\begin{thebibliography}{9}

\bibitem{BBP97}
\textsc{D.H. Bailey}, \textsc{P.B. Borwein} and \textsc{S. Plouffe},
On the rapid computation of various polylogarithmic constants,
\emph{Math. Comput.} \textbf{66} (1997), no. 218, 903--913.

\bibitem{Sa08}
\textsc{V.\,Kh. Salikhov},
On the  irrationality measure of the number $\pi$,
\emph{Russian Math. Surveys} \textbf{63} (2008), no.~3, 570--572.

\bibitem{Sa10}
\textsc{V.\,Kh. Salikhov},
On the measure of irrationality of the number $\pi$,
\emph{Math. Notes} \textbf{88} (2010), no.~4, 563--573.

\bibitem{ZZ20}
\textsc{D. Zeilberger} and \textsc{W. Zudilin},
The irrationality measure of $\pi$ is at most $7.103205334137\dots$,
\emph{Exp. Math.} 24 \textbf{4} (2015), no.~4, 419--423.

\end{thebibliography}
\end{document}